\documentstyle[12pt]{article} 
\newtheorem{theorem}{Theorem} 
\newtheorem{lemma}[theorem]{Lemma} 
 
\newtheorem{proposition}[theorem]{Proposition} 
\def\beq{\begin{equation}}  \def\eeq{\end{equation}} 
\def\bb{\begin{eqnarray*}}  \def\ee{\end{eqnarray*}}
\def\b{\begin{eqnarray}}    \def\e{\end{eqnarray}}
\def\1{\hbox{\rm\setbox1=\hbox{1}\copy1\kern-.5\wd1 I}} 
\def\D{{\hbox{\rm\setbox1=\hbox{I}\copy1\kern-.45\wd1 D}}}
\def\E{{\hbox{\rm\setbox1=\hbox{I}\copy1\kern-.45\wd1 E}}} 
\def\N{\hbox{\rm\setbox1=\hbox{I}\copy1\kern-.45\wd1 N}}
  \def\n{{\cal N}}     \def\L{{\cal L}}
   \def\nin{\in{\!\!\!\!\!/}\,} 
\def\p{\hbox{\rm\setbox1=\hbox{I}\copy1\kern-.45\wd1 P}} 
\def\R{\hbox{\rm\setbox1=\hbox{I}\copy1\kern-.45\wd1 R}}
\def\st{\stackrel{d}{=}}    \def\d{d_{_{TV}}\!}  
\def\bib{\vspace{-3mm}\bibitem}    
   \def\pr{\noindent{\bf Proof}\ }  \def\Z{{\bf Z}_+}  
\def\z{\zeta}   \def\l{\lambda} \def\vp{\varphi}   \def\ve{\varepsilon}
\def\a{\alpha}  \def\({\left(}  \def\){\right)} 

\begin{document} \date{\small 28.07.2020} 
 \author{S.Y.Novak\\ {}\\ {\small MDX University London}} 
 \title{On the T-test} \maketitle 

 \begin{abstract} 
  The $T$-test is probably the most popular statistical test; it is routinely recommended by the textbooks. 
	
  The applicability of the test relies upon the validity of normal or Student's approximation to the distribution of Student's statistic $\,t_n$. However, the latter assumption is not valid as often as assumed. 

  We show that normal or Student's approximation to $\,\L(t_n)\,$ does not hold uniformly even in the class $\,{\cal P}_n\,$ of samples from zero-mean unit-variance bounded distributions. We present lower bounds to the corresponding error. 

  The fact that a non-parametric test is not applicable uniformly to samples from the class $\,{\cal P}_n\,$ seems to be established for the first time. It means the $T$-test can be misleading, and should not be recommended in its present form. 
	
	We suggest a generalisation of the test that allows for variability of possible limiting/approximating distributions to $\,\L(t_n)$.   \end{abstract} 


				\section{Introduction}

  Testing a hypothesis concerning the mean $\,\E X\,$ of the unknown distribution is one of the major tasks in statistical hypotheses testing. 
  In particular, one can be interested in testing the hypothesis 
	$$ H_0=\{\E X\!=\!a\}\,\ \hbox{vs}\,\ H_A=\{\E X\!=\!b\}, $$ 
where $\,a\!\ne\!b$; the hypothesis $\,H_0=\{\E X\!\le\!a\}\,$ vs the hypothesis $\,H_A=\{\E X\!\ge\!b\}$, the hypothesis $\,H_0=\{\E X\!=\!a\}\,$ vs the hypothesis $\,H_A=\{\E X\!\ne\!a\}$, etc. (see, e.g., Lehmann \cite{L}). 
	
	Suppose that $\,X,X_1,X_2,...,X_n\,$ are independent and identically distributed (i.i.d.) random variables (r.v.s). Set $$\,S_n = X_1\!+\!...\!+\!X_n.$$ 
		Throughout the paper we assume that $\,\E X^2<\infty,$ and denote $\,\sigma^2=\hbox{var}\,X$. Below a bar over a random variable means that it is centered by its mathematical expectation. 
	
  The natural estimator of $\,\E X\,$ is the sample mean $$ \hat X = S_n/n, $$ 
leading to the test with the test-statistic $\,(\hat X\!-\!a)\sqrt{n}/\sigma\,$ if $\,\sigma\,$ is known (the so-called $Z$-test) or to the test with the test-statistic 
  $$\,(\hat X\!-\!a)\sqrt{n}/\hat\sigma_n\,,$$ 
where $\,\hat\sigma_n\,$ is an estimator of the standard deviation (the $T$-test). 

  The $T$-test is arguably the most popular statistical test. In view of the law of large numbers (LLN) and the central limit theorem (CLT) it appears perfectly justified if $\,\E X^2<\infty.$ 
	
	However, we show below that the $T$-test has problems even in the simplest case where the main hypothesis is $\,H_0=\{\E X\!=\!a\},$ the alternative hypothesis is $\,H_A=\{\E X\!=\!b\}\,$ ($a\!<\!b$) and $\,\hbox{var}\,X\,$ is known. We argue that the $T$-test is not automatically applicable, and requires prior checks.\\ 

  Let $\,\Phi\,$ denote the standard normal distribution function (d.f.), 
	$$ \vp=\Phi',\ \,\Phi_c=1\!-\!\Phi.	$$ 	
  Given the type-I error $\,\ve\!\in\!(0;1)$, the critical value $\,c_\ve\,$ (and hence the critical region) of the ``two-sided'' $Z$-test of $\,H_0=\{\E X\!=\!a\}\,$ vs $\,H_A=\{\E X\!\ne\!a\}$ is determined by the equation $\,\ve = 2\Phi_c(c_\ve),$ i.e., 
	\beq \label{T0} 
	c_\ve = \Phi_c^{-1}(\ve/2), 
	\eeq 
the acceptance region	of the ``two-sided'' test is $\,[a\!-\!c_\ve\sigma/\!\sqrt{n};a\!+\!c_\ve\sigma/\!\sqrt{n}\,],$ so that 
 \beq \label{T1}  
 \p_a( \hat X\!\nin[a\!-\!c_\ve\sigma/\!\sqrt{n};a\!+\!c_\ve\sigma/\!\sqrt{n}\,] ) 
 = \p_a( (\hat X\!-\!a)\sqrt{n}/\sigma\!\nin[-c_\ve;c_\ve] ) \approx \ve 
 \eeq 
in line with CLT, where $\,\p_a\!=\!\L(X)\,$ assuming $\,\E X\!=\!a$ 
(in the case of the $T$-test $\,\sigma\,$ is unknown and is to be replaced by $\,\hat\sigma_n$; the tradition suggests using Student's d.f. instead of $\,\Phi$). 

  Here $\,[a\!-\!c_\ve\sigma/\!\sqrt{n};a\!+\!c_\ve\sigma/\!\sqrt{n}\,]\,$ is the asymptotic confidence interval based on CLT. 
	A more robust approach suggests using estimates of the accuracy of normal approximation, i.e., replacing the asymptotic confidence interval 
with a sub-asymptotic confidence interval (see \cite{N11}, ch. 9).  

  In the case of a ``one-sided'' $Z$-test of $\,H_0=\{\E X\!=\!a\}\,$ vs $\,H_A=\{\E X\!=\!b\}$, where say $\,a\!>\!b,$ the critical value $\,c_\ve\,$ is determined by the equation $\,\ve = \Phi_c(c_\ve)\,$ and $\,H_0\,$ is rejected if $\,\hat X\!<\!a\!-\!c_\ve\sigma/\!\sqrt{n}\,,$ so that 
	$$ 
	\p_a((\hat X\!-\!a)\sqrt{n}/\sigma\!<-\!c_\ve) \approx \ve \eqno(\ref{T1}^*)
	$$ 
(in the case of the $T$-test $\,\sigma\,$ is replaced by $\,\hat\sigma_n\,$ and $\,\Phi\,$ is replaced by Student's d.f.). 

  Set $\,\delta=(a\!-\!b)/\sigma,$ and let 
	$$ \z_n = (\hat X\!-\!\E X)\sqrt{n}\,/\sigma = \bar S_n/\sigma\sqrt{n}\,. $$ 
  The probability of the type-II error in the case of (\ref{T1}) equals 
  \beq \label{T2} 
	\p_b(\delta\sqrt{n}-\!c_\ve\!\le\!\z_n\!\le\!\delta\sqrt{n}+\!c_\ve) = 
  \p(\z_n\!\ge\!\delta\sqrt{n}-\!c_\ve) - \p(\z_n\!\ge\!\delta\sqrt{n}+\!c_\ve). 
	\eeq 
  In the case of (\ref{T1}$^*$) the probability of the type-II error equals 
	$$ \p_b( (\hat X\!-\!a)\sqrt{n}/\sigma\!\ge\!-c_\ve ) = 
	\p(\z_n\!\ge\!\delta\sqrt{n}-\!c_\ve). \eqno(\ref{T2}^*) $$ 
  Thus, the probabilities of the type-II errors are the large deviations probabilities. 

	The need to approximate the probability of the type-II error led to the rise of the theory of large deviations (see, e.g., \cite{Lin,P65,P95} and references therein).   In view of (\ref{T2}) one needs to approximate the asymptotics of probabilities $\,\p(\z_n\!\ge\!x_n),$ where $\,\{x_n\}\uparrow\,$ is a sequence of positive numbers, $\,x_n=O(\sqrt{n}\,)\,$ as $\,n\!\to\!\infty$. 
	
	Under certain assumptions on $\,\L(X)\,$ and the rate of $\,x_n\,$ 
	\beq \label{T10} 
	\p(\z_n\!\ge\!x_n) \sim \Phi_c(x_n) \qquad(n\!\to\!\infty). 
	\eeq 
  In more general situations (in particular, if $\,x_n\!\asymp\!\sqrt{n}\,$) the asymptotics of $\,\p(\z_n\!\ge\!x_n)\,$ can be expressed in terms of the so-call ``rate function'' $\,\Lambda\,$ (Legendre transform of function $\,\psi(t) = \ln\E e^{tX}\,$ --- see, e.g., Petrov \cite{P65} or \cite{N11}, ch. 14.5). 
	
	One would prefer normal approximation (\ref{T10}) as the rate function is typically unknown, and the task of estimating it can be demanding. Note that one can use the Erd\"os--R\'enyi maximum of partial sums as an estimator of $\,\Lambda^{-1}\,$ (cf. \cite{N11}, ch. 2); however, the accuracy of estimation is likely to be pour. 
	
  An important result is due to Linnik (\cite{Lin}, Theorem 3). 
	Suppose that $\,X,X_1,X_2,...$ are i.i.d.r.v.s, $\,\E X\!=\!0.$	Let $\,\rho(n)\,$ be a monotone function on $\,(0;\infty)\,$ such that $\,0\!<\!\rho(n)\!\uparrow\!\infty,$ $\,\rho(n)\!\ll\!n^\ve\,$ for any $\,\ve\!>\!0$. If 
  \beq \label{Li}  \E\exp(|X|^{4\a/(2\a+1)}) < \infty \eeq 
for some $\,0\!<\!\a\!<\!1/6,$ then 
  \beq \label{Lin} \sup_{0\le x\le n^\a\!/\rho(n)} 
	\left| \p(\z_n\!\ge\!x)/\Phi_c(x)-1 \right| \to0.
	\eeq 
  Moreover, if (\ref{Lin}) holds with $\,n^\a\!/\rho(n)\,$ replaced by 
$\,n^\a\!\rho(n),$ then (\ref{Li}) holds true \cite{Lin}. 

  Thus, the question appears answered in the case of known $\,\sigma\,$ and ``close alternatives'' (i.e., $\,\delta\!\equiv\!\delta(n)\!\to\!0\,$ and $\,1\!\ll\!x_n\!\ll\!\sqrt{n}\,$ as $\,n\!\to\!\infty$). However, we argue that the use of normal approximation is not properly justified. 
  The reason for that is that the test is effectively applied as a non-parametric one --- textbooks implicitly assume that the $T$-test ``works'' uniformly over the non-parametric class $\,{\cal P}_{\sigma}(a_1,a_2)\,$ of distributions with mean $\,\E X\!\in\![a_1;a_2]\,$ and standard deviation $\,\sigma$. 
	
	We show below that weak convergence of $\,\L(\bar S_n/\sqrt{n})\,$ to the normal law cannot hold uniformly in the class $\,{\cal P}_{1}(0,0)\,$ of zero-mean unit-variance distributions (the issue with uniform convergence is known in the literature though not in the context of the $T$-test --- see \cite{N18} and references therein concerning weak convergence uniformly in a class of distributions). 

  Note that this problem does not arise if one deals with a typical parametric family of distributions since a typical parametric family $\,\{P_\theta,\theta\!\!\in\!\!{\bf\Theta}\}\,$ has one-to-one correspondence between a parameter and a distribution. 

  In applications the standard deviation is usually unknown and has to be replaced by its estimator, e.g., $\,\hat\sigma_n\,$ or $\,\tilde\sigma_n\,$, where 
  $$
	\hat\sigma_n^2 = n^{-1}\!\sum^n_{i=1}\!X_i^2\!-\!\hat X^2,\ 
  \tilde\sigma_n^2 = \sum^n_{i=1}(X_i\!-\!\hat X)^2/(n\!-\!1), 
	$$ 
$\,\z_n\,$ is to be replaced with Student's statistic $\,t_n\equiv t_n(X_1,...,X_n),$ where  
	$$
	g t_n = (\hat X\!-\!\E X)\sqrt{n}/\hat\sigma_n = \bar S_n/\hat\sigma_n \sqrt{n}\,. 
	$$  

  The test of the hypothesis $\,H_0=\{\E X\!=\!a\}\,$ involving test statistic $\,t_n\,$ is called the $T$-test. It is one of the most widely used statistical tests. 
	
	Textbooks advocate using the $T$-test when testing for the hypothesis $\,H_0=\{\E X\!=\!a\}\,$ vs the alternative hypothesis $\,H_A=\{\E X\!=\!b\}$, where $\,a\!\ne\!b$; when testing for the hypothesis $\,H_0=\{\E X\!\le\!a\}\,$ vs the hypothesis $\,H_A=\{\E X\!\ge\!b\}$, etc.. 
	
	The test suggests accepting $\,H_0=\{\E X\!=\!a\}\,$ if $\,\hat X\!\in\![a\!-\!c_\ve\hat\sigma_n/\!\sqrt{n} ; a\!+\!c_\ve\hat\sigma_n/\!\sqrt{n}\,]\,$ (or if $\,\hat X\!\in\![a\!-\!c_\ve\tilde\sigma_n/\!\sqrt{n} ; a\!+\!c_\ve\tilde\sigma_n/\!\sqrt{n}\,]$), where $\,c_\ve\,$ is given as in (\ref{T0}) with Student's d.f. instead of $\,\Phi$, so that the probability of the type-I error is asymptotically $\,\ve\,$ (assuming normal or Student's approximation to $\,\L(t_n)$, cf. (\ref{St})). 
	
	Let where $\,\hat\delta=(a\!-\!b)/\hat\sigma_n\,$. 
	The probability of the type-II error equals 
	\bb 
	&& \p_b(\hat X\!\in\![a\!-\!c_\ve\hat\sigma_n/\!\sqrt{n} ; 
	a\!+\!c_\ve\hat\sigma_n/\!\sqrt{n}\,]) = 	\p(\hat\delta\sqrt{n}-\!c_\ve \le 
	\bar S_n/\hat\sigma_n\sqrt{n} \le \hat\delta\sqrt{n}+\!c_\ve)\\ && 
	\p(t_n\!\ge\!\hat\delta\sqrt{n}-\!c_\ve) - 
	\p(t_n\!\ge\!\hat\delta\sqrt{n}+\!c_\ve) 
	\ee 
in the case of a two-sided test (say, $\,a\!>\!b$) or $\,\p(t_n\!\ge\!\hat\delta\sqrt{n}-\!c_\ve)\,$ in the case of a one-sided test. 
  Note that $\,\p_b(t_n\!\ge\!\hat\delta\sqrt{n}-\!c_\ve) \le \p_b(\hat\sigma_n\!>\!2\sigma) + \p_b(t_n\!\ge\!\delta\sqrt{n}/2-\!c_\ve).$ 
  These are the large deviations probabilities.

Since the test involves $\,\bar S_n,$ we may assume in the sequel that $\,\E X\!=\!0$. 
  Denote 
 $$ T_n = \sum^n_{i=1} X_i^2\,,\ t_n^* = S_n/T_n^{1/2}\,. $$ 
 Self-normalised sum $\,t_n^*\,$ is closely related to Student's statistic $\,t_n$: 
 \beq \label{St} 
 t_n\!=\!t_n^*/\sqrt{1\!-{t_n^*}^2/n}\ ,\ t_n^*\!=\!t_n/\sqrt{1\!+t_n^2/n}\,. 
 \eeq 
Note that 
  $$ 
	\{t_n\!\ge\!x\} = \left\{\,t_n^*\ge x/\sqrt{1\!+\!x^2\!/n}\,\right\}\!,\ 
  \{t_n^*\!\ge\!y\} = \left\{\,t_n\ge y/\sqrt{1\!-\!y^2\!/n}\,\right\} 
	\eqno{(\ref{St}^*)} 
	$$ 
if $\,x\!\ge\!0,\,0\!\le\!y\!\le\!\sqrt{n}\,.$ 
Thus, probabilities of the events involving $\,t_n\,$ that appear in the test can be presented as probabilities of the events involving $\,t_n^*\,$. 
  In particular, the limiting distributions of $\,t_n\,$ and $\,t_n^*\,$ coincide. 
	In the sequel we will mainly speak about $\,t_n^*$. 

	Student's statistic converges weakly to the standard normal law if and only if $\,\L(X)\,$ is in the domain of attraction of a normal law and $\,EX\!=\!0\,$ (Gin\'e et al. \cite{GGM}). 
  The class $\,\L_{\cal S}\,$ of limiting distributions of Student's statistic in the case of a triangle array of i.i.d. in each row r.v.s has been described by Mason \cite{M05}. 
  An estimate of the accuracy of normal approximation to the distribution of Student's statistic with explicit constants has been given by Novak \cite{N00,N04} in the case of identically distributed r.v.s and by Shao \cite{Shao05} in the case of non-identically distributed r.v.s (see also \cite{N11,Pi16} and references therein). 
	
		A (\ref{Lin})-type result for the self-normalised sum $\,t_n^*\,$ has been suggested as well \cite{JSW}: there exists an absolute constant $\,A\,$ such that 
	\beq \label{JSW} 
	\left| \p(t_n^*\!\ge\!x) - \Phi_c(x) \right| \le 
	A(1\!+\!x)^2e^{-x^2/2} \E|X|^3/\sigma^3\sqrt{n}\, 
	\qquad(0\!\le\!x^3\!\le\!\sigma^3\sqrt{n}/\E|X|^3). 
	\eeq  
In view of (\ref{St6}), bound (\ref{JSW}) is equivalent to 
	$$ \left| \p(t_n^*\!\ge\!x)/\Phi_c(x) - 1 \right| \le 
	A(1\!+\!x)^3 \E|X|^3/\sigma^3\sqrt{n}\, \ \ \ \ \ 
	\qquad(0\!\le\!x^3\!\le\!\sigma^3\sqrt{n}/\E|X|^3). \eqno(\ref{JSW}^*) 
	$$ 
	
  As in the case of known $\,\sigma,$ one can easily get an impression that it is save to apply the $T$-test to any sample of i.i.d. observations with a finite second moment. 
	
	Recall that the $T$-test was originally formulated for samples of i.i.d. normal r.v.s. 	In most applications the observations are not normally distributed. 
	Nonetheless, textbooks suggest that normal approximation is applicable if the sample size is large: ``the size of the one- and two-sample $T$-tests is relatively insensitive to nonnormality (at least for large samples). Power values of the $T$-tests obtained under normality are asymptotically valid also for all other distributions with finite variance. This is a useful result...'' (\cite{L}, p. 207). 
	
	This opinion appears widely accepted suggesting that the probabilities of the type-I and type-II errors can be accurately approximated using normal or Student's d.f.. 

	The intuition behind such a suggestion is obvious: 
	$$ t_n^* = \frac{S_n/\sigma\sqrt{n}}{T_n/n\sigma^2}\,. $$ 	
  By CLT, $\,S_n/\sigma\sqrt{n}\,$ is expected to converge weakly to a standard normal r.v., while by LLN $\,T_n/n\sigma^2\,$ is expected to converge to 1 as $\,n\!\to\!\infty\,$ (cf. \cite{L}, p. 205). 
	
	The purpose of this article is to show that such an impression can be misleading. What textbooks are missing is that weak convergence of $\,\L(t_n)\,$ and $\,\L(t_n^*)\,$ to the normal law is not uniform. 
  Theorems \ref{Th1}, \ref{Th2} show that the $T$-test cannot be applied uniformly even in the class of bounded zero-mean unit-variance distributions.\\ 
		
  Let $\,\Psi_n\,$ denote the distribution function of Student's statistic with $\,n\,$ degrees of freedom. The use of Student's distribution instead of normal one has been inherited from the case of normally distributed observations. However, it is easy to check that $\,\Psi_n\,$ is close to $\,\Phi$: 
	\beq \label{*} 
	\sup_x|\Psi_n(x)-\Phi(x)| \le C/n \qquad(n\!\to\!\infty)
	\eeq
(cf. Pinelis \cite{Pi13}). 
  The table of Student's distribution function shows little difference between $\,\Psi_n(\cdot)\,$ and $\,\Phi(\cdot)\,$ if $\,n\!\ge\!60$. 
	Thus, preference to $\,\Psi_n\,$ over $\,\Phi\,$ appears questionable. 
	
	We argue that normal or Student's approximation to the distribution of Student's statistic is not automatically applicable. 
	We suggest performing prior checks in order to find out if a particular (not necessarily normal) approximation to the distribution of the test statistic is applicable. This leads to a generalisation of the $T$-test that allows for non-conventional approximating distributions. 
	We discuss implications for the choice of critical levels. 

  Section 2 addresses the question of validity of the $T$-test uniformly over a class of distributions with a finite variance. 
	Section 3 presents an example of non-normal approximation as well as an estimate of the accuracy of such approximation in terms of the total variation distance. 
	The approximating distribution appears new in the literature on the topic. 
  Section 4 suggests a generalisation of the $T$-test. Proofs are postponed to section 5.

				\section{Problems with the $T$-test}

  The $T$-test has been criticized by a number of authors. 
For instance, Bahadur (\cite{Bah71}, Example 8.1) shows that Student's $T$-test is not Bahadur-efficient if $\,H_0=\{\E X\!=\!0\}\,$ and $\,X_1,...,X_n\,$ are i.i.d. normal $\,\n(\theta;1)\,$ r.v., $\,\theta\!\ge\!0$. 
	
  Rukhin \cite{Ruh} shows that Student's $T$-test is not Bahadur-efficient in the case of testing the null hypothesis $\,H_0=\{\theta\!=\!0\}\,$ against $\,H_1=\{\theta\!=\theta_1\!\}\,$ for a location-scale parameter family $\,F_{\theta,c}$, where $\,F_{\theta,c}(x) = F((x\!-\!\theta)/c)\ (\forall x),$ $\,\theta\!\in\!\R,\,c\!>\!0,$ $\,F\,$ is a d.f. with a finite (in a neighbourhood of 0) moment generating function. 

  The aim of this article is to show that normal or Student's approximation to the distribution of Student's statistic is not automatically applicable, and the test can be misleading. 
	
	Textbooks effectively suggest applying the $T$-test as a non-parametric test; the class of distributions considered applicable is effectively the class of all distributions with finite variances. 
  We show below that the use of the $T$-test is not justified in such generality 
even in the case of testing a simple hypothesis $\,H_0=\{\E X\!=\!a\}\,$ against a simple alternative $\,H_A=\{\E X\!=\!b\}\,$ in the assumption that $\,\hbox{var}\,X\!<\!\infty$. 	

	W.l.o.g. we may assume in the sequel that $\,a\!=\!0.$ 

  Let $\,{\cal P}_n\,$ denote the class of distributions $\,{\cal L}(X_1,...,X_n)\,$ of random vectors $\,(X_1,...,X_n)\,$ such that $\,X,X_1,...,X_n\,$ are independent and identically distributed bounded random variables, $\,\E X\!=\!0,$ $\,\E X^2\!=\!1$. 
  The use of normal approximation in the $T$-test would be justified if normal approximation held uniformly in the class $\,{\cal P}_n$. 

  We show below that normal approximation is not applicable uniformly in the class $\,{\cal P}_n.$ In particular, there exists an absolute constant $\,c\!>\!0\,$ such that for any $\,n\!>\!3$ 
	\beq \label{T8} 
	\inf_{x\ge0} \sup_{{\cal P}_n} \left| \p(t_n^*\!\ge\!x)/\Phi_c(x)-1 \right| 
	\ge c. 
	\eeq 
	A similar result holds if $\,\Phi\,$ in (\ref{T8}) is replaced with $\,\Psi_n\,$ or $\,\Psi_{n-k}$, where $\,k\!\in\!\N$. 	
		
	A comparison of (\ref{JSW}$^*$) with Linnik's result (\ref{Lin}) suggests $\,\L(t_n^*)\,$ has ``better'' asymptotic properties than $\,\L(\z_n)$. However, it has been noticed in \cite{N04} that $\,\L(t_n^*)\,$ has certain disadvantages comparing to $\,\L(\z_n)$.   In particular, a non-uniform Berry-Esseen-type inequality does not hold for Student's statistic (though a modified non-uniform Berry-Esseen-type inequality is valid, see Theorem 12.24 in \cite{N11}).

	Example 12.3 in \cite{N11} presents a situation where 
  \beq \label{T9} 
	\sup_{{\cal P}_n} \left| \p(t_n^*\!\ge\!x_n)/\Phi_c(x_n)-1 \right| \to\infty
	\eeq 
as $\,n\!\to\!\infty\,$ if $\,x_n\!=\!\sqrt{n}\,.$ 
  	An natural question is if 
		$$ \sup_{{\cal P}_n} \left| \p(t_n^*\!\ge\!x_n)/\Phi_c(x_n)-1 \right|\to0 $$ 
as $\,n\!\to\!\infty\,$ for a particular sequence $\,\{x_n\}\,$ such that $\,0\!\le\!x_n\!\ll\!\sqrt{n}\,$. 

  Theorems \ref{Th1}, \ref{Th2} below answer that question. 
In particular, we show that the $T$-test is not applicable uniformly over $\,{\cal P}_n\,$ regardless of the size of the sample. In other words, the outcome of the test can be misleading even for large-size samples.

  \begin{theorem} \label{Th1} For any $\,n\!>\!3$ 
	\beq                                                   \label{T3} 
  \inf_{x\ge0}\sup_{{\cal P}_n} \Big| \p(t_n^*\!\ge\!x)/\Phi_c(x) - 1 \Big| \ge 
	1.25 e^{-1/2(n\!-\!2)}-1>0. 
  \eeq 
  If $\,\{x_n\}\,$ is a non-decreasing sequence of positive numbers such that $\,1\!\ll\!x_n\!\le\!\sqrt{n}\,$ as $\,n\!\to\!\infty,$ then 
  \beq                                                   \label{T3+}
	\sup_{{\cal P}_n} \Big| \p(t_n^*\!\ge\!x_n)/\Phi_c(x_n) - 1 \Big| \to\infty 
	\qquad(n\!\to\!\infty). \eeq 
	\end{theorem} 

	A similar result holds if normal approximation to $\,\L(t_n^*)\,$ has been replaced with Student's approximation. 
  Denote 	$$ \Psi_n^c = 1\!-\!\Psi_n\,,\ \psi_n=\Psi_n'. $$

 \begin{theorem} \label{Th2} 
As $\,n\!\to\!\infty,$ 
	\beq                                                   \label{T6}
  \inf_{x\ge0}\sup_{{\cal P}_n} \Big| \p(t_n^*\!\ge\!x)/\Psi_n^c(x) - 1 \Big| 
	\ge 1/4\!+\!o(1). 
	\eeq 
If $\,\{x_n\}\,$ is a non-decreasing sequence of positive numbers such that $\,1\!\ll\!x_n\!\ll\!\sqrt{n}\,$ as $\,n\!\to\!\infty,$ then 
  $$ \sup_{{\cal P}_n} \Big| \p(t_n^*\!\ge\!x_n)/\Psi_n^c(x_n)-1 \Big| 
	\to\infty \qquad(n\!\to\!\infty). \eqno(\ref{T6}^*) $$ 
	\end{theorem}

  The result holds if $\,\Psi_n\,$ in (\ref{T6}) has been replaced with $\,\Psi_{n-k}\,,$ where $\,k\,$ is a fixed natural number.

			\section{An example of non-normal approximation}

  It may be counter-intuitive to expect that Poisson (or Binomial) distribution may play any role in the study of the properties of the $T$-test but Proposition \ref{Th3} below states it may.

	In this section we present an example of non-normal/non-Student's approximation to $\,\L(t_n)\,$ and $\,\L(t_n^*)\,$ and evaluate the accuracy of the approximation. 
	The example highlights the fact that the limiting distribution of Student's statistic may take on value $\,\infty\,$ with positive probability. 

  Given r.v.s $\,Y\,$ and $\,Z$, we denote by $\,\d(Y;Z)\equiv\d(\L(Y);\L(Z))\,$ the total variation distance between $\,\L(Y)\,$ and $\,\L(Z)$. 

  Let $\,\pi_\l\,$ denote a Poisson r.v. with parameter $\,\l$. Set 
	\beq \label{T19} 
	Y = (np-\pi_{np})\Big/\!\sqrt{np^2\!+\!(1\!-\!2p)\pi_{np}}\,. 
	\eeq 
Note that $$\,\p(Y\!=\!\sqrt{n}) = e^{-np}\,.$$ 
	Given $\,p\!\in\!(0;1],$ we set $\,q=1\!-\!p,$ 
$\,\ve_n\!=\!\min\!\left\{1; \(2\pi[(n\!-\!1)p]\)^{-1/2}\! + 
2(1\!-\!e^{-np})p/(1\!-\!1/n)\right\}.$

  \begin{proposition} \label{Th3} 
	Let $\,X,X_1,...,X_n\,$ be i.i.d.r.v.s with the distribution 
 \beq                                         \label{Ex12.3}
 \p\(X\!=\!\sqrt{p/q}\,\) = q,\ \p\(X\!=\!-\sqrt{q/p}\,\) = p, 
 \eeq 
where $\,p\!\in\!(0;1/2].$ Then 
 \beq                                                   \label{Ex123} 
 \d\(t^*_n;Y\) \le 3p/4e + 2(1\!-\!e^{-np})^2p^2 + 2(1\!-\!e^{-np})p^2\ve_n.
 \eeq 
  \end{proposition} 

  In the light of (\ref{St}), inequality (\ref{Ex123}) can be reformulated as follows: 
  $$ \d\(t_n;\eta\) \le 3p/4e + 2(1\!-\!e^{-np})^2p^2 + 2(1\!-\!e^{-np})p^2\ve_n, 
	\eqno(\ref{Ex123}^+) $$ 
where $\,\eta = (np\!-\!\pi_{np})\Big/\!\sqrt{\pi_{np}(1\!-\!\pi_{np}/n)}\,.$

  Denote 
  $$ Y_\l = (\l\!-\!\pi_\l)/\sqrt{\pi_\l}\,\qquad(\l\!>\!0). $$ 
Clearly, $\,Y_\l\,$ is a defective random variable: $Y_\l\,$ takes on value $\,\infty\,$ with probability $\,e^{-\l}$. According to Proposition \ref{Th3}, 
	\beq \label{T11} t_n \Rightarrow Y_\l\,,\ t^*_n \Rightarrow Y_\l	\eeq 
if $\,p=p(n)\!\sim\!\l/n\,$ as $\,n\!\to\!\infty.$ 

	Weak convergence (\ref{T11}) may hold in more general situations, e.g., if $\,X_i\!\st\!(\xi_i\!-\!\E\xi)/\E^{1/2}\xi\,$ and $\,\xi_1,\xi_2,...,\xi_n\,$ are i.i.d. non-degenerate r.v.s taking values in $\,\Z\!=\!\{0,1,2,...\}.$ 
	For example, (\ref{T11}) holds if $\,X_i\st (p\!-\!\eta_i)/\sqrt{p}\,,$ where $\,\{\eta_i\}\,$ are i.i.d. Poisson $\,{\bf\Pi}(p)\,$ r.v.s with $\,p=p(n)\!\sim\!\l/n\,$ as $\,n\!\to\!\infty.$ 
	
	In situations where $\,\L(t^*_n)\,$ can be approximated by $\,\L(Y)\,$ or $\,\L(Y_\l),$ the ``asymptotic approach'' suggests the critical values $\,c_-,c_+\,$ can be chosen according to the equations 
	$$ 
	\p(Y\!>\!c_+) = \p(Y\!<\!c_-) = \ve/2\ \ \hbox{or}\ \ 
	\p(Y_\l\!>\!c_+) = \p(Y_\l\!<\!c_-) = \ve/2 
	$$ 
with $\,\l\!=\!np\,$ replaced by its consistent estimator (the ``two-sided'' test); the ``sub-asymptotic approach'' suggests incorporating estimate (\ref{Ex123}). 

  A possible alternative to distribution $\,\L(X)\,$ given by (\ref{Ex12.3}) is $\,\L(X') = \L(X\!+b),$ where $\,b\!\ne\!0.$ If $\,(X_1',...,X_n')_{n\ge1}\,$ is a triangular array of i.i.d.r.v.s such that $\,X_i'\st b\!+\!(p\!-\!\xi_i)/\sqrt{pq}\,,$ then the probability $\,\p(c_-\!\le\!t_n^*(X_1',...,X_n')\!\le\!c_+)\,$ of the type-II error is a probability of large deviations for the Binomial distribution.

								\section{A generalised test}

  The $T$-test relies on the validity of normal (or Student's) approximation to $\,\L(t_n)$.  
	The common impression is that $\,\L(t_n)\,$ is close to the standard normal distribution if the sample size $\,n\,$ is large (see, e.g., Lehman \cite{L}, p. 205). 

	It is known that the limiting distribution of $\,t_n\,$ is not always normal (see Mason \cite{M05}). 

  In this section we suggest a generalised $T$-test. The idea is to check first if a particular approximation (not necessarily normal or Student's) is applicable. 
	The latter can be done using sharp estimates of the accuracy of approximation. 

  Thus, the generalised $T$-test requires (1) a list of possible limiting/approximating distributions; (2) sharp estimates of the accuracy of approximation of $\,\L(t_n)\,$ by the corresponding distributions; 
	(3) estimation of certain quantities involved in those estimates of the accuracy of approximation (e.g., estimation of $\,\sigma\,$ and $\,\E|X^3|\,$ in the case of normal approximation). 

  Traditionally, the obvious candidate for the approximating distribution is the standard normal law. 
	One can employ the following approximate bound to the uniform distance between $\,\L(t_n^*)\,$ and $\,\n(0;1)$ (cf. \cite{N04}, Corollary 2): 
  \beq \label{Cor2} 
  |\p\(t_n^*\!<\!x\)-\Phi(x)| \le (6.4\hat\mu_3/\hat\sigma^3 + 
  2\hat\mu_1/\hat\sigma)/\sqrt{n}\,, 
  \eeq 
where $\,\hat\mu_k\,$ denotes a consistent estimator of $\,\mu_k := \E|X\!-\!\E X|^k,$ $\,k\!\ge\!1$; we denote $\,\hat\sigma^2\!:=\!\hat\mu_2$. 

  Bound (\ref{Cor2}) is based on the estimate of the accuracy of normal approximation to $\,\L(t_n^*)\,$ from \cite{N04} that seems to be the sharpest available in the case of i.i.d. observations (cf. the discussion in \cite{Pi16}, Remarks 4.16--4.17). 
	
	The bound in \cite{N04} involves a term (say, $\,\gamma_n\,$) which is of order $\,o(n^{-1/2})$. In applications moments $\,\{\mu_k\}\,$ have to be substituted by their consistent estimators, generating an extra error. Therefore, it is reasonable to omit the term $\,\gamma_n\,$ and arrive at (\ref{Cor2}). 
	
  The use of normal approximation can be justified if the right-hand side (r.h.s.) of (\ref{Cor2}) is less than a certain small number (say, $\,\ve$) specified by a statistician (e.g., $\,\ve\!=\!0.01$). 

  Since the limiting distribution of $\,t_n\,$ may differ from $\,\n(0;1)\,$ (cf. Proposition \ref{Th3}), we suggest that one first checks if a particular (not necessarily normal) approximation to the distribution of the test statistic $\,t_n\,$ is applicable. 
	
	One may have a number of bounds of the type 
	\beq \sup_x |\p(t_n^*\!\le\!x) - F_k(x)| \le r_n(k), \eeq 
where $\,\{F_k\}\,$	are d.f.s of certain candidate distributions. It is natural to choose $\,k\!=\!k_*\,$ such that $\,r_n(k_*)=\min_k r_n(k)$. 

  Obviously, one needs a list of possible limiting/approximating distributions together with the corresponding estimates of the accuracy of approximation with explicit constants. Such a list will always be finite but until recently only normal and Student's distributions were on the list. 	
	
	Proposition \ref{Th3} adds another candidate to the list. 

  Note that one can have a situation where neither distribution from the list has the estimate of the accuracy of approximation (e.g., $\,r_n(k)$) below the specified threshold level $\,\ve\,$ (i.e., $\,\min_k r_n(k)>\ve$). That would mean the $T$-test is not applicable (either because of a small sample size or because of the list been too short).

				\section{Proofs}

  Since $\,t_n\,$ and $\,t_n^*\,$ are scale-invariant, w.l.o.g. we may assume in the sequel that $\,\hbox{var}\,X=1.$ 
		Below the operation of multiplication is superior to the division. 

  The proofs of Theorems \ref{Th1}, \ref{Th2} use the fact that $\,\L(t_n)\,$ and $\,\L(t_n^*)\,$ are not stochastically bounded uniformly in $\,{\cal P}_n.$\\ 

  \pr of Theorem \ref{Th1}. Note that $\,t_n^*\!\le\!\sqrt{n}\,.$ Thus, (\ref{T3}) trivially holds if $\,x_n\!>\!\sqrt{n}\,$. 
	Therefore, we may assume that $\,x\!\in\![0;\sqrt{n}\,]$. 
	
	It suffices to find i.i.d. bounded r.v.s $\,X,X_1,...,X_n\,$ such that $\,\E X\!=0,$ $\,\E X^2\!=1,$ and (\ref{T3}) holds. 
  We will employ distribution (\ref{Ex12.3}) that seems to play the role of a testing stone when one deals with self-normalised sums and Student's statistic. 

  Let $\,X\,$ be a r.v. with distribution (\ref{Ex12.3}). Then 
	$$ X_i \st (p-\xi_i)/\sqrt{pq}\,, \eqno(\ref{Ex12.3}^*) $$ 
where $\,\{\xi_i\}\,$ are independent Bernoulli $\,{\bf B}(p)\,$ r.v.s (cf. Example 12.3 in \cite{N11}). Note that 
$$ \E X=0,\ \E X^2=1,\ \E |X|^3 = (p^2\!+\!q^2)/\sqrt{pq}\,. $$ 
In particular, $\,{\cal L}(X_1,...,X_n)\in{\cal P}_n.$ 

  Denote $\,S_n^\xi = \xi_1+...+\xi_n\,.$ Then  
	\beq \label{T15} 
	S_n = (np\!-\!S_n^\xi)/\sqrt{pq}\,,\ T_n = np/q+(1\!-\!2p)S_n^\xi/pq. 
	\eeq 
Hence	$$ t_n^* = (np-S_n^\xi)/\sqrt{np^2\!+\!(q\!-\!p)S_n^\xi}\ . $$ 
	
	Set 
	\beq \label{T21} g(k) = (np\!-\!k)/\sqrt{np^2\!+\!(q\!-\!p)k}\,\qquad(k\!\in\!\Z). 
	\eeq
Note that $\,t_n^* = g(S_n^\xi).$ 
  Since function $\,g(\cdot)\downarrow\,,$ we have 
	\beq \label{Extr} \p(t_n^*\ge g(k)) = \p(S_n^\xi\le k). \eeq 
 
  Clearly, $\,t_n^*\,$ takes on its largest possible value $\,g(0)=\sqrt{n}\,$ when $\,X_1=...=X_n= \sqrt{p/q}\,,$ $\,t_n^*\,$ takes on its second largest possible value $\,g(1) = (np\!-\!1)/\sqrt{np^2\!+\!q\!-\!p}\,$ when $\,n\!-\!1\,$ sample elements equal $\,\sqrt{p/q}\,$ and one sample element equals $\,-\sqrt{q/p}\,,$ etc.. Hence 
 \beq  \label{St6+} 
 \p(t_n^*\!=\!\sqrt{n}\,) = q^n,\ 
 \p\Big(t_n^*\!=\!(np\!-\!1)/\sqrt{np^2\!+\!(q\!-\!p)}\,\Big) = npq^{n-1}. \eeq

  We consider first the case where $\,x\!\in\![0;1]$. 
According to (\ref{Extr}), (\ref{St6+}), 
  $$ \p(t_n^*\ge g(1)) = (q\!+\!np)q^{n-1}. $$ 
  Note that 
	$$ \ln(1\!-\!x) \ge -x\!-\!x^2/2(1\!-\!x)^2\qquad(0\!\le\!x\!<\!1). $$ 
Hence \beq \label{T7} (1\!-\!p)^n \ge \exp(-np(1\!+\!p/2q^2)). \eeq 
	
  Denote 
  $$ p_x = \(1\!+x\sqrt{1\!-\!1/n}\Big/\sqrt{1\!-\!x^2/n}\,\)\!\!\Big/n\,. $$ 
Set $\,p=p_x\,.$ Then $\,g(1)=x.$ 

  One can check that $\,np/q \ge 1\!+\!x.$ Hence 
  $$ \p(t_n^*\ge x) \ge (2\!+\!x)q^n. $$ 
Taking into account (\ref{T7}), we derive 
  \bb
  \p(t_n^*\ge x) &\ge& (2\!+\!x) 
	\exp\!\(-\Big(1\!+\!x\sqrt{1\!-\!1/n}\Big/\sqrt{1\!-\!x^2/n}\,\Big)(1\!+\!p/2q^2)\)\\ 
	&\ge& (2\!+\!x) \exp\!\Big(-(1\!+\!x) \Big( 1\!+\!(1\!+\!x)/2n(1\!-\!2/n)^2 \Big)\Big)
	. \ee 
	
  It is well-known that $\,\Phi_c(x)\le e^{-x^2/2}\!/2.$ Hence 
  \bb 
  \p(t_n^*\ge x)/\Phi_c(x) &\ge& 2(2\!+\!x) 
	\exp(x^2\!/2\!-\!1\!-\!x\!-\!(1\!+\!x)^2/2n(1\!-\!2/n)^2)\\ 
  &\ge& \frac{_2}{^e} (2\!+\!x) \exp(x^2\!/2\!-\!x\!-\!2/n(1\!-\!2/n)^2. 
  \ee 
  Note that function $\,h(x) = x^2/2\!-\!x\!+\ln(2\!+\!x)\,$ takes on its minimum in $\,[0;1]\,$ at $\,x_*=(\sqrt{5}-\!1)/2 \approx 0.618.$ Hence 
$\,\frac{_2}{^e} (2\!+\!x) \exp(x^2/2\!-\!x) > 1.256.$ Thus, 
  $$ 
	\p(t_n^*\ge x)/\Phi_c(x) > 1.25 e^{-2/n(1\!-\!2/n)^2} \qquad(n\!>\!3). 
	\eqno(\ref{T3}^*) 
	$$ 
In particular, $\,\p(t_n^*\ge x)/\Phi_c(x)>1.01\,$ if $\,n\!>\!12$.\\ 

  We consider now the case where $\,x\!\in\![1;\sqrt{n}\,]$. It is well-known that 
 \beq                                               \label{St6} 
 \frac1{1\!+\!x} < 
 \frac {\Phi_c(x)}{\vp(x)} < \frac 1x \qquad(x\!>\!0). 
 \eeq
Relations (\ref{St6+}) -- (\ref{St6}) yield 
  $$
	\p(t_n^*\!\ge\!x)/\Phi_c(x) \ge \p(t_n^*\!\ge\!\sqrt{n})/\Phi_c(x) \ge 
	(1\!-\!p)^n x/\vp(x). $$ 
	
	Let $\,p=1/n.$ Then 
  \beq                                               \label{T12} 
  \p(t_n^*\!\ge\!x)/\Phi_c(x_n) \ge 
	\frac{_{\sqrt{2\pi}\,}}{^e} x e^{x^2/2-1/2(n\!-\!2)}\,. 
  \eeq 
Since $\,\inf_{x\ge1} xe^{x^2/2} = e^{1/2}\,,$ we have 
 $$ 
 \p(t_n^*\!\ge\!x)/\Phi_c(x) \ge 
 \frac{_{\sqrt{2\pi}\,}}{^{\sqrt{e}}} e^{-1/2(n\!-\!2)}\,. 
 $$ 
Note that $\,\sqrt{2\pi/e}>1.52$. 
Thus, (\ref{T3}) holds. In particular,  
  $$ 
  \sup_{{\cal P}_n} (\p(t_n^*\!\ge\!x)/\Phi_c(x)\!-\!1) 
	\ge \sqrt{2\pi}/e^{3/4}-1 \eqno(\ref{T3}^\star)  
	$$ 
if $\,n\!>\!3$. 
Relation (\ref{T3+}) follows from (\ref{T12}). 
  The proof is complete. \hspace*{\fill} $\Box$\\

  {\it Remark 1}. The statement of Theorem \ref{Th1} can be reformulated for negative $\,x\,$ by switching from $\,\{X_i\}\,$ to $\,\{-X_i\}$: (\ref{T3}) implies that for any $\,n\!>\!3$ 
	$$  
  \inf_{x\le0}\sup_{{\cal P}_n} \Big| \p(t_n^*\!\le\!x)/\Phi(x) - 1 \Big| \ge 1.25 e^{-2/n(1\!-\!2/n)^2}. \eqno(\ref{T3}^+) 
  $$ 
Similarly one reformulates the statement of Theorem \ref{Th2}:
as $\,n\!\to\!\infty,$ 
	$$    
  \inf_{x\le0}\sup_{{\cal P}_n} \Big| \p(t_n^*\!\le\!x)/\Psi_n(x) - 1 \Big| 
	\ge 1/4\!+\!o(1). \eqno(\ref{T6}^*) $$  

   {\it Remark 2}. Distribution (\ref{Ex12.3}) is not the only one that can be used in order to establish (\ref{T8}). 
	For instance, let independent r.v.s $\,\tau\,$ and $\,\eta\,$ be independent of $\,\xi,$ $\,\L(\tau)={\bf B}(c/n),$ where $\,c\!\ge\!0,$ 
	$\,\E\eta=0,$ $\,\E\eta^2=1.$ Set 
		$$ X = \tau\eta + (1\!-\!\tau)(p-\xi)/\sqrt{pq}\,, $$ 
and let $\,\{X_i\}\,$ be independent copies of $\,X$. 
  Then $\,\E X=0,$ $\,\E X^2=1.$ 
	
Let, for example, $\,x=0.$ If $\,p=1/n,$ then 
  $$ 
	\p(t_n^*\ge0)/\Phi_c(0) \ge (1\!-\!c/n)^n q^{n-1}(q\!+\!np)	\sim 2/e^{1+c} 
	$$ 
as $\,n\!\to\!\infty$. Therefore, $\,\p(t_n^*\!\ge\!0)/\Phi_c(0) \ge 4/e^{1+c}\!+\!o(1) \!>\!1\,$ for all large enough $\,n\,$ if $\,c<\ln(4/e)$.\\

  \pr of Theorem \ref{Th2} involves Lemma \ref{Th4} and the argument from the proof of Theorem \ref{Th1}. 
	
Since $\,t_n^*\!\le\!\sqrt{n}\,,$ (\ref{T6}) trivially holds if $\,x_n\!>\!\sqrt{n}\,$. Below we may assume that $\,x\!\in\![0;\sqrt{n}\,]$. 
  Let $\,X,X_1,...,X_n\,$ be as in the proof of Theorem \ref{Th1}. Recall that 
  $$ \psi_n(x) = C_n (1\!+\!x^2/n)^{-(n+1)/2}\qquad(x\!\in\!\R), $$ 
where 
	\beq \label{T13} C_n = \Gamma((n\!+\!1)/2)/\sqrt{\pi n}\,\Gamma(n/2),\ 
	\Gamma(y)\!=\!\int_0^\infty\! t^{y-1}e^{-t}dt\qquad(y\!>\!0). \eeq 
	
	We consider first the case where $\,x\!\in\![1;\sqrt{n}\,]$. 
		Using (\ref{T5}), we derive 
	\b	\nonumber 
	\p(t_n^*\!\ge\!x)/\Psi_n^c(x) &\ge& 
	\p(t_n^*\!\ge\!\sqrt{n})/\Psi_n^c(x)\\ \label{T14} 
	&\ge& (1\!-\!1/n)^{n+1} x(1\!+\!x^2/n)^{(n-1)/2}\!\big/C_n  
	\e 
if $\,p=1/n.$ 
  It is known that $\,C_n\!\to\!1/\sqrt{2\pi}\,$ as $\,n\!\to\!\infty$. Since $\,\inf_{x\ge1} x(1\!+\!x^2/n)^{(n-1)/2} = (1\!+\!1/n)^{(n-1)/2} = \sqrt{e}+\!o(1),$ we have 
	$$ 
	\p(t_n^*\!\ge\!x)/\Psi_n^c(x) \ge \sqrt{2\pi/e}+\!o(1) \qquad(n\!\to\!\infty)
	$$ 
uniformly in $\,x\!\in\![1;\sqrt{n}\,]$. 

  We consider now the case where $\,x\!\in\![0;1]$. Let $\,p=p_x\,.$ 
Then 
$$ 
\p(t_n^*\ge x) \ge (2\!+\!x) e^{-1-x}(1\!+\!o(1)) \qquad(n\!\to\!\infty) 
$$ 
uniformly in $\,x\!\in\![0;1]$.   Taking into account (\ref{T20}), it is easy to see that $\,\Psi_n^c(x)\to\Phi_c(x)\,$ as $\,n\!\to\!\infty\,$ uniformly in $\,x\!\in\![0;1]$. Therefore,  
  $$ 
	\p(t_n^*\!\ge\!x)/\Psi_n^c(x) \ge (2\!+\!x) e^{-1-x}/\Phi_c(x) 
	(1\!+\!o(1)) \qquad(n\!\to\!\infty) $$ 
uniformly in $\,x\!\in\![0;1]$. 
  Repeating the argument that led to (\ref{T3}$^*$), we derive 
$$ \p(t_n^*\!\ge\!x)/\Psi_n^c(x) > 1.25+o(1) \qquad(n\!\to\!\infty). 
$$ 
Thus, $\,\inf_{x\ge0}\sup_{{\cal P}_n} \Big| \p(t_n^*\!\ge\!x)/\Psi_n^c(x) - 1 \Big| \ge 1/4\!+\!o(1).$ 

  If $\,\{x_n\}\,$ is a non-decreasing sequence of positive numbers such that $\,1\!\ll\!x_n\!\ll\!\sqrt{n}\,$ as $\,n\!\to\!\infty,$ then (\ref{T14}) entails (\ref{T6}$^*$). The proof is complete. \hspace*{\fill}$\Box$\\

 \begin{lemma} \label{Th4} As $\,n\!>\!1,x\!>\!0,$ 
 \beq                                                   \label{T5} 
 \frac{\sqrt{2\pi}\,C_n}{\sqrt{1\!+\!1/n}} \Phi_c\Big(x\sqrt{1\!+\!1/n}\,\Big) \le 
 \Psi_n^c(x) \le \frac{C_n}{(1\!-\!1/n)x} (1\!+\!x^2/n)^{-(n-1)/2}\ .\\ 
 \eeq \end{lemma}

   Note that (\ref{T5}) means $\,\Psi_n^c(x)\,$ decays rather fast when $\,x\!\in\!(0;\sqrt{n}\,]$: 
  $$ 
	\!\!\!\!\!\Psi_n^c(x) \le C_n e^{-\frac{\,x^2}4(1\!-\!1/n)}\!/x(1\!-\!1/n), 
	\eqno(\ref{T5}') $$ 
	$$ \Psi_n^c(x) \ge C_n e^{-\frac{\,x^2}2(1+1/n)}\!/(1\!+\!x)(1\!+\!1/n). 	\eqno(\ref{T5}'') 
	$$  
	
  \pr of Lemma \ref{Th4}. It is easy to see that 
 \bb 
 \Psi_n^c(x) \!&=&\! C_n\int_x^\infty (1\!+\!y^2/n)^{-(n+1)/2}\, dy\\  
 \!&\le&\! C_n x^{-1}\! \int_x^\infty (1\!+\!y^2/n)^{-(n+1)/2}\, ydy\\ 
 \!&=&\! \frac{C_n}{(1\!-\!1/n)}\, x^{-1} (1\!+\!x^2/n)^{-(n-1)/2}\ . 
 \ee 
  Using Taylor's formula, one can check that 
	\beq \label{T20} y\ge \ln(1\!+\!y) \ge y-y^2/2 \qquad(y\!\ge\!0). \eeq 
Hence 
  $$ 
	e^{x^2} \ge (1\!+\!x^2/n)^{n} \ge \exp(x^2-x^4/2n) \ge e^{x^2/2} 
	\qquad(0\!\le\!x^2\!\le\!n). 
	$$ 
Therefore, 
  $$ 
	\Psi_n^c(x) \le \frac{C_n}{(1\!-\!1/n)} x^{-1} e^{-x^2(1-1/n)/4}\,. 
	$$ 
Similarly, 
 \bb 
 \Psi_n^c(x) \!&\ge& C_n \int_x^\infty \exp(-y^2(1\!+\!1/n)/2) dy\\ 
 \!&=&\! C_n \sqrt{2\pi/(1\!+\!1/n)}\, \Phi_c\Big(x\sqrt{1\!+\!1/n}\,\Big)\\ 
 \!&\ge& C_n e^{-x^2(1+1/n)/2} \!/(1\!+\!x)(1\!+\!1/n) 
 \ee 
by (\ref{St6}). The proof is complete. \hspace*{\fill} $\Box$\\

  \pr of Proposition \ref{Th3}. 
	Recall that r.v.s $\,\{X_i\}\,$ obey (\ref{Ex12.3}$^*$) and 
  $$ t_n^* = (np-S_n^\xi)/\sqrt{np^2\!+\!(q\!-\!p)S_n^\xi}\ , $$ 
where $\,S_n^\xi = \sum_{i=1}^n \xi_i\,.$ Note that 
  \beq \label{T4} t_n^* = g(S_n^\xi),\ Y = g(\pi_{np}), \eeq 
where monotone function $\,g\,$ is given by (\ref{T21}). 

  Theorem 4.12 in \cite{N11} states that 
	\beq                                     	\label{RN} 
	\d(S_n^\xi;\pi_{np}) \le 3p/4e + 2\delta^2 + 2\delta^*\ve_n, 
  \eeq 
where $$\,\delta = (1\!-\!e^{-np})p,\ \delta^* = (1\!-\!e^{-np})p^2. $$ 
	For any $\,A\!\subset\!\Z\,$ let $\,B:=g(A).$ Taking into account (\ref{T4}), we observe that 
	$$ 
	\p(t^*_n\!\in\!A)-\p(Y_{n,p}\!\in\!A) = 
\p(g(S_n^\xi)\!\in\!B)-\p(g(\pi_{np})\!\in\!B) \le \d(S_n^\xi;\pi_{np}). 
	$$  
Thus, (\ref{Ex123}) follows from (\ref{RN}). The proof is complete. \hspace*{\fill} $\Box$\\ 

	{\bf Conclusion}. 
	We have shown that the $T$-test in its present form can be misleading even if the sample size is arbitrarily large: normal or Student's approximation to the distribution of Student's statistics $\,t_n\,$ is not automatically applicable. 
	
	Note that the sample size is always finite; in applications it often cannot be increased either due to physical restrictions or because of cost considerations. 

  The paper suggests a generalisation of the $T$-test that involves checking for the appropriate approximating distribution, and requires estimates with explicit constants of the accuracy of approximation to $\,\L(t_n)$. 
	
	The list of possible limiting/approximating distributions may include, beyond normal, functions of Poisson, compound Poisson, and some other infinitely divisible laws (cf. (\ref{T19})). 
	
	The problem of deriving estimates of the accuracy of normal approximation 
with explicit constants to the distribution of a sum of r.v.s goes back to Tchebychef \cite{Cheb86} and Liapunov \cite{Lia}. 
	It lead to a vast literature with contributions from many renowned authors (see, e.g., references in \cite{AZ,N11,P95,Shev13}). 
  The task of evaluating the accuracy of Poisson and compound Poisson approximation has been addressed by many distinguished authors (see, e.g., references in \cite{AZ,C16,Pre83,Zai03}). 
	
  In 1950s Kolmogorov has formulated the problem of evaluating the accuracy of approximation of the distribution of a sum of independent r.v.s by infinitely divisible laws. 
	The topic has attracted a lot of attention among researchers (see, e.g., references in \cite{AZ,C16,M05,Zai03}). However, most estimates obtained so far have implicit constants. In particular, the task of evaluating the constants in Arak's, Presman's and Zaitsev's inequalities remains open.  
	
	We conclude that the generalised $T$-test requires a list of possible limiting/approximating distributions together with the corresponding estimates of the accuracy of approximation. 
	The class $\,\L_{\cal S}\,$ of limiting distributions of Student's statistic has been described by Mason \cite{M05}. 
	For most distributions from $\,\L_{\cal S}\,$ the task of deriving estimates of the accuracy of approximation with explicit constants remains open. 

				\end{document}